\documentclass[12pt,english]{article}
\usepackage{theorem}
\usepackage{babel}
\usepackage{mathrsfs}
\usepackage{a4wide}
\usepackage{graphics}
\usepackage{epsfig}
\usepackage{amssymb}
\usepackage{latexsym}
\usepackage{amsmath}

\newtheorem{thh}{Theorem}[section]
\newtheorem{df}[thh]{Definition}

\newtheorem{lem}[thh]{Lemma}
\newtheorem{cor}[thh]{Corollary}

\title{Invariant hypersurfaces for derivations \\ in positive characteristic}
\author{Philippe Bonnet}
\date{}

\newcommand{\dem}{{\em Proof: }}
\newcommand{\qed}{\begin{flushright} $\blacksquare$\end{flushright}}

\newcommand{\DD}{{\cal{D}}}
\newcommand{\FF}{{\cal{F}}}
\newcommand{\XX}{{\cal{X}}}
\newcommand{\ZP}{\mathbb{Z}/p}
\newcommand{\FP}{\mathbb{F}_p}
\newcommand{\GAF}{\Pi(A,{\cal{D}})}

\begin{document}
\maketitle

\begin{center}
Mathematisches Institut, Universit\"at Basel \\
Rheinsprung 21, 4051 Basel, Switzerland \\ e-mail:
Philippe.bonnet@unibas.ch
\end{center}

\begin{abstract} 
Let $A$ be an integral $k$-algebra of finite type over an algebraically closed field $k$ of characteristic $p>0$. Given a collection $\DD$ of
$k$-derivations on $A$, that we interpret as algebraic vector fields on $X=Spec(A)$, we study the group spanned by the hypersurfaces
$V(f)$ of $X$ invariant under $\DD$ modulo the rational first integrals of $\DD$. We prove that this group is always a finite dimensional
$\FP$-vector space, and we give an estimate for its dimension. This is to be related to the results
of Jouanolou and others on the number of hypersurfaces invariant under a foliation of codimension 1.
As a application, given a $k$-algebra $B$ between $A^p$ and $A$, we show that the kernel of the
pull-back morphism $Pic(B)\rightarrow Pic(A)$ is a finite $\FP$-vector space. In particular, if
$A$ is a UFD, then the Picard group of $B$ is finite.
\end{abstract}

\section{Introduction}

Let $A$ be an integral $k$-algebra of finite type over an algebraically closed field $k$ of characteristic $p>0$, and let $K(A)$ be
its fraction field. Denote by $X$ the spectrum of $A$, i.e. $X=Spec(A)$. Let $Der_k(A)$ be the set of $k$-derivations of $A$, or in
other words the set of algebraic vector fields on $X$.
Given a collection $\DD$ of $k$-derivations on $A$, we would like to compare the set of hypersurfaces invariant under $\DD$ with its set of rational
first integrals, i.e. the rational functions on $X$ annihilated by $\DD$.

We begin with a few definitions and recalls about the theory of foliations in positive characteristic (see \cite{M-P} or \cite{Ek}). The
space $Der_k(A)$ is provided with a Lie bracket defined by the rule $[D_1,D_2]=D_1 D_2 - D_2 D_1$. Moreover if $D$
is a $k$-derivation on $A$, then $D^p=D\circ ...\circ D$ is again a $k$-derivation on $A$. This follows easily from
the Leibniz formula with binomial coefficients for $D^p(fg)$, since $char(k)=p$. {\em A foliation $\FF$ on $X$} is a sub-$A$-module of
$Der_k(A)$, stable by Lie bracket and $p$-closed, that is:
\begin{itemize}
\item{for any $D,D'$ in $\FF$, $[D,D']=D\circ D' - D'\circ D$ belongs to $\FF$,}
\item{for any $D$ in $\FF$, $D^p=D\circ ...\circ D$ belongs to $\FF$.}
\end{itemize}
The definition of a foliation in characteristic zero is exactly the same without the $p$-closedness condition. The codimension of $\FF$ is
defined as $codim(\FF)=dim(A)- rk(\FF)$, where $rk(\FF)$ is the rank of $\FF$ as an $A$-module. For convenience,
given a hypersurface $H$ of $X$, denote by $I_H$ the ideal of elements of $A$ which vanish along $H$.

\begin{df}
Let $\DD$ be a subset of $Der_k(A)$. The hypersurface $H$ is invariant by $\DD$
if the ideal $I_H$ is stable by all elements of $\DD$, i.e. $D(I_H)\subseteq I_H$
for all $D\in \DD$. The element $f$ of $K(A)^*$ is a rational first integral if
$D(f)=0$ for all $D\in \DD$.
\end{df}
Geometrically speaking, this means that all algebraic vector fields $D$ in $\DD$ are tangent to $H$.
If the ground field is $\mathbb{C}$, Jouanolou proved that, under some conditions,
a codimension 1 foliation has either finitely many invariant hypersurfaces, or has a nonconstant first integral (see \cite{Jou}).
This result was latter improved by Brunella and Nicolau for complex codimension 1 foliations and for $A$-submodules
$\DD$ of $Der_k(A)$ of corank 1 in positive characteristic (see \cite{B-N}). More precisely, they proved that such an
$A$-submodule with infinitely many irreducible invariant hypersurfaces must have a non trivial first integral $f$, i.e.
$f\not\in K(A)^p$. The proof consists in showing first that $\DD$ is a codimension 1 foliation, and then use the fact that
foliations have always non trivial first integrals in positive characteristic (see \cite{M-P}).

However this result does not give any information
about the relationship between invariant hypersurfaces and first integrals.
Our purpose is to study this relationship when we restrict to hypersurfaces of the form $V(f)$, where $f$ belongs to $A$. By analogy with the
theory of foliations (see \cite{Jou2}), we give the following:

\begin{df}
Let $\DD$ be any collection of $k$-derivations on $A$. An element $f$ of $K(A)^*$ is an algebraic solution of $\DD$ if $D(f)/f$ belongs to $A$
for any $D$ in $\DD$.
\end{df}
This definition extends the notion of invariant hypersurface in the following sense: if $A$ is a UFD and $f\in A$ is irreducible,
then $f$ is an algebraic solution if and only if the ideal $(f)$ is stable by $\DD$, or in other words if
$V(f)$ is invariant by $\DD$. Note that any first integral is also an algebraic solution. Since:
$$
\frac{D(fg)}{fg}=\frac{D(f)}{f} +\frac{D(g)}{g}
$$
for any $f,g\in K(A)^*$, the set of algebraic solutions (resp. first integrals) forms a multiplicative group. In order to compare
both notions, we introduce the quotient:
$$
\Pi(A,\DD)=\frac{ \{\mbox{Algebraic solutions of }\DD \}}
{\{\mbox{First integrals of }\DD\}}
$$
Since $\Pi(A,\DD)$ is an abelian group, it is a $\mathbb{Z}$-module. But $char(k)=p$, so $D(f^p)=0$
for any $f\in K(A)$ and any $D\in Der_k(A)$. Therefore $\GAF$ is a $\FP$-vector space, where
$\FP=\ZP$. We are going to prove the following:

\begin{thh} \label{solutions}
The $\FP$-vector space $\GAF$ is finite dimensional.
\end{thh}
In particular, the space of algebraic solutions differs from the first integrals by only a finite set. This theorem is to be related
with a result of Pereira (see \cite{Pe}) which asserts that, for $A=k[x_1,...,x_n]$, every general algebraic vector field $D$ has
a nontrivial algebraic solution $f$, i.e $f\in A-A^p$. In a sense, a general algebraic vector field in $k[x_1,...,x_n]$
has nontrivial algebraic solutions, but not too many compared to its first integrals. 

Using the arguments of the proof of Theorem \ref{solutions}, one can also derive an estimate for the dimension of $\GAF$
if $\DD$ is a foliation.
Let $\deg$ be a degree function on $K(A)$ (see Section \ref{degrees}). Given a $k$-derivation $D$ on $A$, we define the
{\em degree $\deg(D)$} of $D$ as: 
$$
\deg(D)=\sup_{f \in K(A)}\{\deg(D(f))-\deg(f)\}
$$
For instance, if $A=k[x_1,...,x_n]$ and $\deg$ stands for the standard homogeneous degree on $A$, then the degree of a $k$-derivation $D$ is
given by $\deg(D)=\sup\{\deg(D(x_i))\} - 1$. In this case, the degree $\deg(D)$ of any $k$-derivation $D$ is finite. The boundedness of the degree
for a $k$-derivation is analogous to the notion of continuity of a derivation with respect to a valuation (see \cite{Mo}).

\begin{df}
Let $\deg:K(A)\rightarrow \mathbb{Z}\cup\{-\infty\}$ be a degree such that $deg(a)=0$ for any element $a\in k^*$. Then $\deg$
is good if the following conditions hold:
\begin{itemize}
\item{For any $k$-derivation $D$ on $A$, the degree of $D$ is finite.}
\item{For any integer, the dimension $l(n)$ of $\{f\in A\; | \; \deg(f)\leq n\}$ over $k$ is finite.}
\end{itemize}
\end{df}

\begin{cor} \label{estimate}
Let $A$ be an integral $k$-algebra of finite type over an algebraically closed field $k$ of characteristic $p$. Let $\deg$ be a good degree function on $K(A)$. Let
$\FF$ be a foliation on $Spec(A)$, spanned by $D_1,...,D_n$. Then we have:
$$
dim_{\FP} \Pi(A,\FF) \leq l(\deg(D_1))+...+l(\deg(D_n))
$$
\end{cor}
We derive from Theorem \ref{solutions} another result that has a priori nothing to do with foliations.
Given two integral $k$-algebras $A,B$ such that $A^p\subseteq B\subseteq A$, we would like to find a relationship between their Picard groups. Denote
by $\FF_B$ the maximal foliation on $Spec(A)$ which vanishes on $B$, i.e. the set of $k$-derivations on $A$ which vanish
on $B$. If $A$ and $B$ are normal rings, then $B$ is exactly the kernel of $\FF_B$ (see \cite{M-P}). Using this
interpretation, we can prove the following:

\begin{cor} \label{picard}
Let $A$ be a normal integral $k$-algebra of finite type over an algebraically closed field $k$ of characteristic $p$. Let $B$ be a
normal subalgebra of $A$ such that $A^p \subseteq B \subseteq A$. Then the kernel of the
pull-back morphism $\pi: Pic (B) \rightarrow Pic(A), \; M\mapsto M\otimes A$ is a finite $\FP$-vector space, and we have:
$$
dim_{\FP} \ker\pi \leq dim_{\FP} \Pi(A,\FF_B)
$$
In particular, if $A$ is a UFD, then the Picard group of $B$ is finite.
\end{cor} 
We end up this paper with two examples. The first one illustrates the sharpness of the estimate given in Corollary \ref{estimate},
and the necessity to consider foliations in its formulation. We will use the second one to determine a Picard group.

\section{Reduction to the case of foliations}

Let $A$ be an integral $k$-algebra of finite type, where $k$ is algebraically closed of characteristic $p$.
In this section, we are going to see how to restrict ourselves to finite collections of derivations satisfying some further conditions.
For any collection $\DD$ of $Der_k(A)$, denote by $\FF(\DD)$ {\em the smallest foliation containing} $\DD$, i.e. the intersection of all
foliations on $X$ containing $\DD$. By definition, it is well-defined and unique.

\begin{lem} \label{reduction}
Let $A$ be an integral $k$-algebra of finite type, where $char(k)=p$, and let $\DD$ be a collection of $k$-derivations on $A$.
Let $\{D_1,...,D_n\}$ be a system of generators of $\FF(\DD)$. Then $\DD$ and $\{D_1,...,D_n\}$ have exactly the same set of algebraic
solutions (resp. first integrals). In particular $\GAF=\Pi(A,\{D_1,...,D_n\})$.
\end{lem}
\dem For any collection $\DD_0$, denote by $AS(\DD_0)$ the group of algebraic solutions of $\DD_0$, and by $FI(\DD_0)$ its group of first
integrals. Our purpose is to show the equalities:
$$
AS(\DD)=AS(\{D_1,...,D_n\})\quad \mbox{and} \quad FI(\DD)=FI(\{D_1,...,D_n\})
$$
We only need to prove the first one, since the proof for the second is entirely similar.
So let $f$ be an algebraic solution for $\{D_1,...,D_n\}$. Then every element $D$ of $\DD$ is an $A$-linear combination
of $D_1,...,D_n$. Since $D_i(f)/f$ belongs to $A$ for any $i$, $D(f)/f$ belongs to $A$. Since this holds for any $D\in \DD$,
$f$ is an algebraic solution of $\DD$ and we have:
$$
AS(\{D_1,...,D_n\})\subseteq AS(\DD)
$$
To prove the equality of both sets, it suffices to show that every algebraic solution of $\DD$ is also an algebraic solution
of $\{D_1,...,D_n\}$. So fix any algebraic solution $f$ of $\DD$. Consider the set $M(f)$ of $k$-derivations on $A$ such that $D(f)/f$
belongs to $A$. By definition $M(f)$ is an $A$-submodule of $Der_k(A)$ containing $\DD$, and $f$ is an algebraic solution of $M(f)$. We
prove that $M(f)$ is stable by Lie bracket and $p$-closed. First, let $D_1,D_2$ be any elements of $M(f)$. Then there exist some elements
$R_{1},R_{2}$ of $A$ such that $D_1(f)=R_{1}f$ and $D_2(f)=R_{2}f$, and we have:
$$
[D_1,D_2](f)=\{D_1(R_{2}) - D_2(R_{1})\} f
$$
So $[D_1,D_2]$ belongs to $M(f)$. Second, let $D$ be any element of $M(f)$. Then there exists an element $R$ of $A$ such that $D(f)=Rf$.
Consider the sequence $\{R_{n}\}_{n>0}$ in $A$ constructed by induction as follows:
$$
R_{1}=R \quad \mbox{and} \quad \forall n>0, \; R_{n+1}= D(R_{n}) + RR_n
$$
By induction, we have $D^n(f)=R_n f$ for any $n>0$. In particular $D^p(f)/f$ belongs to $A$ and $D^p$ belongs to $M(f)$. Therefore $M(f)$
is a foliation containing $\DD$. By definition of $\FF(\DD)$, $M$ contains $\FF(\DD)$, hence its generators $D_1,...,D_n$. Since $f$ is
an algebraic solution of $M(f)$, $f$ is also an algebraic solution of $\{D_1,...,D_n\}$, and the result follows.
\qed

\section{A lemma on degree functions} \label{degrees}

In this section, we are going to establish a lemma that is crucial for the proof of theorem \ref{solutions}. This lemma
asserts that every integral $k$-algebra $A$ of finite type carries a finite set of degree functions enjoying some nice
properties. In fact, this result holds for any field $k$, but for convenience we will only prove it when $k$ is algebraically
closed. In what follows, all degree functions we consider are maps $\deg:K(A)\rightarrow \mathbb{Z} \cup \{-\infty\}$
satisfying the usual axioms, with the additional condition that $\deg(f)=0$ for any $f\in k^*$.

\begin{df}
Let $A$ be an integral $k$-algebra, $K(A)$ its fraction field and $\deg$ a degree on $K(A)$. A $k$-derivation $D$ on $A$
is bounded for $\deg$ if $\deg(D)$ is finite.
\end{df}

\begin{lem} \label{degree}
Let $A$ be an integral $k$-algebra of finite type over an algebraically closed field $k$. Then there exist a finite set $\{\deg_1,...,\deg_r\}$
of degree functions on $K(A)$ satisfying the following conditions:
\begin{itemize}
\item{for any $i=1,...,r$, every $k$-derivation $D$ on $A$ is bounded for $\deg_i$,}
\item{for any integers $n_1,...,n_r$, ${\cal{L}}(n_1,...,n_r)=\{f\in A|\forall i=1,...,r, \; \deg_i(f)\leq n_i\}$ is
a $k$-vector space of finite dimension.}
\end{itemize}
\end{lem}
\dem Let $X$ be the affine variety $Spec(A)$. This variety is embedded in some $k^n$. Let $X'$ be the projective closure of $X$
in $\mathbb{P}^n(k)$, and denote by ${\cal{X}}$ its normalization. By construction, the variety ${\cal{X}}$ is projective, normal
and birational to $X$. Let $H'$ be the hyperplane at infinity in $\mathbb{P}^n(k)$ such that $X=X'-H'$. If
$\varphi:{\cal{X}}\rightarrow X'$ denotes the normalisation morphism, set $H=\varphi^{-1}(H')$. Then $H$
is a finite union of prime Weyl divisors $Z_1,...,Z_r$. Let
$ord_{Z_i}$ be the order along $Z_i$. Since ${\cal{X}}$ is birational to $X$, every element $f$ of $K(A)$ can be considered
as a rational function on ${\cal{X}}$. We set:
$$
\deg_i(f)=-ord_{Z_i}(f)
$$
Since $ord_{Z_i}$ is a valuation, $\deg_i$ defines a degree function on $K(A)$. We are going to show
that these degrees enjoy all the conditions of the lemma. \\

\noindent
{\underline{\em First step}}: Given any $k$-derivation $D$ on $A$, we are going to prove that $D$ is bounded for $deg_i$. Consider
an affine open set $U_i$ in ${\cal{X}}$ such that $Z_i \cap U_i\not=\emptyset$. Let $x_1,...,x_n$ be a set of generators
of ${\cal{O}}_{U_i}$. Since the fraction field of ${\cal{O}}_{U_i}$ is equal to $K(A)$, there exist some elements $\alpha_1,...,
\alpha_n$ of $K({\cal{O}}_{U_i})$ such that $D(x_j)=\alpha_j$. Write $\alpha_j=a_j/b_j$, where $a_j,b_j$ belong to ${\cal{O}}_{U_i}$, and
set $B=b_1...b_n$. By construction, the $k$-derivation $BD$ maps ${\cal{O}}_{U_i}$ into itself, and in particular we have:
$$
BD({\cal{O}}_{U_i, Z_i\cap U_i})\subseteq  {\cal{O}}_{U_i, Z_i \cap U_i}
$$
Since $U_i$ is normal, ${\cal{O}}_{U_i, Z_i \cap U_i}$ is a discrete
valuation ring (see for instance \cite{Ha}). Let $h_i$ be a generator of the unique maximal ideal of ${\cal{O}}_{U_i, Z_i \cap U_i}$, and set $\delta_i=ord_{Z_i}(B)$.
By construction, any rational map $f$ on $X$ with $\deg_i(f)=r$ can be written as:
$$
f= \frac{g}{h_i ^r}
$$
where $g$ is invertible in ${\cal{O}}_{U_i, Z_i \cap U_i}$. By derivation and multiplication, we get:
$$
h_i ^{r+1}BD(f)= h_i BD(g) -rBD(h_i)g 
$$
So $h_i ^{r+1}BD(f)$ belongs to ${\cal{O}}_{U_i, Z_i \cap U_i}$, its order along $Z_i$ is nonnegative
and we find:
$$
\deg_i(f) +1 - \delta_i - \deg_i(D(f))\geq 0
$$
Since this holds for any $f\in K(A)$, we obtain that $\deg_i(D)\leq 1 - \delta_i$, hence $\deg_i(D)$ is finite. \\

\noindent
{\underline{\em Second step}}: We are going to prove the second assertion of the lemma. Fix some integers $n_1,...,n_r$. By definition, the hyperplane
divisor $H$ is a linear combination of the $Z_i$ with positive coefficients. So there exists a positive integer $n$ such that $nH\geq n_1Z_1+...+n_rZ_r$.
Let $f$ be any element of ${\cal{L}}(n_1,...,n_r)$, viewed as a regular function on ${\cal{X}}-H$. If $div(f)$ denotes its Weyl divisor on ${\cal{X}}$,
we have the relation:
$$
div(f) + nH \geq div(f) + n_1 Z_1 +...+n_r Z_r\geq 0
$$
In particular, if $k({\cal{X}})$ denotes the fraction field of ${\cal{X}}$, we have the inclusion:
$$
{\cal{L}}(n_1,...,n_r)\subseteq\{f\in k({\cal{X}})|\; div(f)+ nH\geq 0\}
$$
Since $H$ is an hyperplane section, it is locally principal and it defines an invertible sheaf ${\cal{L}}$ on ${\cal{X}}$. Moreover, the locally principal
divisor $div(f)+nH$ on ${\cal{X}}$ is effective if and only if $div(f)+ nH\geq 0$. Indeed since $\XX$ is normal, a rational function is regular
on an open set $U$ of ${\cal{X}}$ if and only if it has no hypersurface of poles on $U$. So the right-hand side of the latter inclusion
corresponds to the space $\Gamma(\XX, {\cal{L}}^{-n})$ of global sections. Since $\XX$ is projective, this space is finite
dimensional by Serre's Theorem (see \cite{Sh}). 
\qed

\section{A lemma on polynomials}

In this section, we are going to establish a lemma on the shape of some polynomials of low degree having roots in $\FP ^s$. This result
can be stated as follows:

\begin{lem} \label{solution3}
Let $A$ be an integral $k$-algebra where $k$ is a field of characteristic $p$. Let $P(t_1,...,t_s)$ be an element of $A[t_1,...,t_s]$,
of degree $\leq p$ in $t_1,...,t_s$, such that $P(z_1,...,z_s)=0$ for any $(z_1,...,z_s)\in \FP ^s$. Then there exist
some unique elements $a_1,...,a_s$ of $A$ such that:
$$
P(t_1,...,t_s)= a_1 (t_1 ^p - t_1) + ...+ a_s (t_s ^p - t_s)
$$ 
\end{lem}
\dem by induction on $s\geq 1$, the case $s=1$ being clear. Indeed if $P(t)$ vanishes
at any element
$z$ of $\FP$, then $P(t)$ is divisible by $(t^p -t)$. Since $P$ has degree $\leq p$ in $t$, $P(t)=a(t^p -t)$
for some element $a\in A$, and this element is unique. So assume the property holds to the order $s-1$, and let
$P(t_1,...,t_s)$ be an element of $A[t_1,...,t_s]$, of degree $\leq p$ in $t_1,...,t_s$, such that
$P(z_1,...,z_s)=0$ for any $(z_1,...,z_s)\in \FP ^s$. Denote by $a_s$ the coefficient of $t_s ^p$ in the
expression of $P$, and set:
$$
Q(t_1,...,t_s)=P(t_1,...,t_s) -a_s(t_s ^p - t_s)
$$
By definition, the polynomial $Q$ has degree $\leq p$ in $t_1,...,t_s$ and vanishes at any point of $\FP ^s$.
Moreover its degree with respect to the variable $t_s$ is $<p$. Therefore it has an expansion of the form:
$$
Q(t_1,...,t_s) = Q_0(t_1,...,t_{s-1}) + ...+ Q_{p-1}(t_1,...,t_{s-1})t_s ^{p-1}
$$
where each polynomial $Q_i$ has degree $\leq p-i$ in $t_1,...,t_{s-1}$. For any fixed element $z=(z_1,...,z_{s-1})$
of $\FP ^{s-1}$, consider the polynomial:
$$
Q_z(t_s)=Q(z_1,...,z_{s-1},t_s)
$$
By construction $Q_z$ vanishes at any point $z_s$ of $\FP$, and has degree $<p$ in $t_s$. So it is of the form
$a(t_s ^p -t_s)$, where $a$ belongs to $A$. Since $Q_z$ has degree $<p$, it is the zero polynomial. In particular,
for any $z=(z_1,...,z_{s-1})$ of $\FP ^{s-1}$, we have:
$$
Q_0(z_1,...,z_{s-1})=Q_1(z_1,...,z_{s-1})=...=Q_{p-1}(z_1,...,z_{s-1})=0
$$
For any index $i\geq 0$, the polynomial $Q_i$ vanishes on $\FP^{s-1}$. Since its degree in $t_1,...,t_{s-1}$
is $\leq (p-i)$, the polynomial $Q_i$ is by induction an $A$-linear combination of $(t_1 ^p -t_1),...,(t_{s-1} ^p -t_{s-1})$.
But for any $i>0$, $Q_i$ has degree $\leq p-i<p$, so it is the zero polynomial. In particular, we find:
$$
Q(t_1,...,t_{s})=Q_0(t_1,...,t_{s-1})=a_1 (t_1 ^p - t_1) + ...+ a_{s-1} (t_{s-1} ^p - t_{s-1})
$$
for some elements $a_1,...,a_{s-1}$ of $A$. By construction of $Q$, this yields:
$$
P(t_1,...,t_s)=a_1 (t_1 ^p - t_1) + ...+ a_s (t_{s} ^p - t_{s})
$$
The uniqueness of $a_1,...,a_s$ is obvious.
\qed

\section{Proof of Theorem \ref{solutions} and Corollary \ref{estimate}}

Let $A$ be an integral $k$-algebra of finite type, where $k$ is algebraically closed of characteristic $p$. Let $\DD$ be a
collection of $k$-derivations on $A$. We fix a finite set $\{D_1,...,D_n\}$ of generators for the foliation $\FF(\DD)$.
By lemma \ref{reduction}, we know that $\Pi(A,\DD)=\Pi(A,\{D_1,...,D_n\})$, and that the following group morphism:
$$
L: \Pi(A,\DD)\longrightarrow A^n\quad , \quad f\longmapsto \left(\frac{D_1(f)}{f},...,\frac{D_n(f)}{f} \right)
$$
is well-defined and injective. In particular, its image $I$ is isomorphic to $\Pi(A,\DD)$. In order to establish Theorem
\ref{solutions} and Corollary \ref{estimate}, we only need to show that $I$ is a finite dimensional $\FP$-vector space, and
estimate its dimension. This is exactly what we will do in the following subsections.

\subsection{Dimension of $Vect_k(I)$}

Let $Vect_k(I)$ be the vector space over $k$ spanned by $I$ in $A^n$. We prove that {\em $Vect_k(I)$ is finite dimensional}, as follows.
Let $\{\deg_1,...,\deg_r\}$ be a collection of degrees on $K(A)$ satisfying the conditions of lemma \ref{degree}. Then there
exist some constant $\{d_{i,j}\}$ such that, for any $f\in K(A)$:
$$
\deg_i(D_j(f))\leq d_{i,j} + \deg_{i}(f)
$$
Let $f$ be any algebraic solution of $\FF$, and set $L_j(f)=D_j(f)/f$ for any $j$. Then $L_j(f)$ belongs to $A$ for any $j$, and
we have $\deg_i(L_j(f))\leq d_{i,j}$ for any $i,j$. In particular, we obtain the inclusion:
$$
I\subseteq \oplus_{j=1} ^{n} {\cal{L}}(d_{1,j},...,d_{r,j})
$$
By lemma \ref{degree}, all the components of this sum are finite-dimensional over $k$. Therefore $I$ spans a finite dimensional
$k$-vector space.

\subsection{Associated polynomials}

By the previous subsection, there exist some algebraic solutions $f_1,...,f_s$ such that the elements $L(f_1),...,L(f_s)$ form a basis
of $Vect_k(I)$. {\em We fix these algebraic solutions from now on}. By definition, given an algebraic solution $f$, the element $L(f)$ is a
$k$-linear combination of $L(f_1),..., L(f_s)$. The problem is, we do not know a priori what are the possible coefficients of this linear combination.
In this subsection, we are going to construct some polynomials in $s$ variables with the following remarkable property: {\em these polynomials
vanish at every point $(z_1,...,z_s)$ such that $z_1 L(f_1) + ... + z_s L(f_s)$ is equal to some $L(f)$, where $f$ is an algebraic solution}.
The construction proceeds as follows. Since $D_1,...,D_n$ span a foliation, there exist some elements $a_{i,j}$ of
$A$ such that, for any $i=1,...,n$:
$$
D_i ^p= \sum_{j=1} ^n a_{i,j} D_j
$$
We fix these elements $a_{i,j}$, and extend the action of $D_i$ to the ring $A[t_1,...,t_s]$ by setting $D_i(t_{\alpha})=0$ for any
$\alpha\in \{1,...,s\}$. For any $g\in A[t_1,...,t_s]$ and any $i=1,...,n$, we introduce the $k$-linear operator:
$$
T_{i,g}: A[t_1,...,t_s] \longrightarrow A[t_1,...,t_s]\quad , \quad f \longmapsto D_i(f) + gf
$$
For convenience, we set $L_i(f)=D_i(f)/f$. Note that $L(f)=(L_1(f),...,L_n(f))$. For any $i=1,...,n$, consider the following element
$P_i$ of $A[t_1,...,t_s]$: 
$$
P_i(t_1,...,t_s)= \left(T_{i,\sum_{\alpha=1} ^s t_{\alpha} L_i(f_{\alpha})}\right) ^{p-1}\left(\sum_{\alpha=1} ^s t_{\alpha} L_i(f_{\alpha})\right) -
\sum_{j=1} ^n a_{i,j} \left(\sum_{\alpha=1} ^s t_{\alpha} L_j(f_{\alpha})\right)
$$
\begin{df}
$P_1,...,P_n$ are the associated polynomials to $\{D_1,..,D_n;f_1,..,f_s\}$.
\end{df}

\begin{lem} \label{solution0}
Let $f$ be a nonzero element of $K(A)$. If we set $g=D_i(f)/f$, then we have $D_i ^p(f)= (T_{i,g})^{p-1}(g) f$.
\end{lem}
\dem It suffices to prove that $D_i ^N(f)= (T_{i,g})^{N-1}(g) f$ for any $N>0$. We do it by induction on $N>0$, the case
$N=1$ being clear by construction. Assume the assertion holds to the order $N-1$. Then we have to the order $N$:
$$
D_i ^N(f)=D_i(D_i ^{N-1}(f))=D_i((T_{i,g})^{N-2}(g) f)=D_i((T_{i,g})^{N-2}(g))f +(T_{i,g})^{N-2}(g) D_i(f)
$$
By definition of $g$, we have $D_i(f)=gf$ and this implies:
$$
D_i ^N(f)=D_i((T_{i,g})^{N-2}(g))f +(T_{i,g})^{N-2}(g)g f=(T_{i,g}\circ (T_{i,g})^{N-2}(g)))f=(T_{i,g})^{N-1}(g)f
$$
\qed

\begin{lem} \label{solution1}
Let $f$ be any algebraic solution of $\{D_1,...,D_n\}$. Then for any $i=1,...,n$, $L_1(f),...,L_n(f)$ satisfy the relation
$(T_{i,L_i(f)}) ^{p-1}(L_i(f))= \sum_{j=1} ^n a_{i,j} L_j(f)$.
\end{lem}
\dem Let $f$ be an algebraic solution of $D_1,...,D_n$. If $g=D_i(f)/f=L_i(f)$, then we find by applying Lemma \ref{solution0} to $f$:
$$
D_i ^p(f)=(T_{i,L_i(f)}) ^{p-1}(L_i(f)) f
$$
Since $D_i ^p= \sum_j a_{i,j} D_j$ for any $i=1,...,n$ by construction, this implies:
$$
\frac{D_i ^p(f)}{f}=(T_{i,L_i(f)}) ^{p-1}(L_i(f))=\sum_{j=1} ^n a_{i,j} \frac{D_j(f)}{f}=\sum_{j=1} ^n a_{i,j} L_j(f)
$$
\qed

\begin{lem} \label{solution2}
Let $f$ be any algebraic solution of $\{D_1,...,D_n\}$. Let $z_1,...,z_s$ be some elements of $k$ such that $L(f)=z_1 L(f_1)+...+z_s L(f_s)$.
Then for any $i=1,...,n$, we have $P_i(z_1,...,z_s)=0$. Moreover $P_1,...,P_n$ vanish at any point $(z_1,...,z_s)$ of
$\FP ^s$.
\end{lem}
\dem Let $f$ be an algebraic solution of $\{D_1,...,D_n\}$. If $L(f)=z_1 L(f_1)+...+z_s L(f_s)$, then we have $L_i(f)=z_1 L_i(f_1)+...+z_s L_i(f_s)$
for any index $i$. By definition of the associated polynomials and by Lemma \ref{solution1}, we get:
$$
\begin{array}{rcl}
P_i(z_1,...,z_s)&= &(T_{i,\sum_{\alpha=1} ^s z_{\alpha} L_i(f_{\alpha}) }) ^{p-1}\left(\sum_{\alpha=1} ^s z_{\alpha} L_i(f_{\alpha})\right)
- \sum_{j=1} ^n a_{i,j} \left(\sum_{\alpha=1} ^s z_{\alpha} L_j(f_{\alpha})\right) \\ \\
&=& (T_{i,L_i(f)}) ^{p-1}(L_i(f))- \sum_{j=1} ^n a_{i,j} L_j(f)=0
\end{array}
$$
and the first assertion follows. For the second, let $z_1,...,z_s$ be some elements of $\mathbb{F}_p$. Let $y_1,...,y_s$ be some integers
whose classes modulo $p$ are equal to $z_1,...,z_s$ respectively. Then the element $f=f_1 ^{y_1} ... f_s ^{y_s}$ is also an algebraic
solution of $\{D_1,...,D_n\}$, and $L(f)=z_1 L(f_1) +...+z_s L(f_s)$. Therefore $P_i(z_1,...,z_s)=0$
for any $i=1,...,n$ by the first assertion.
\qed

\begin{lem} \label{solution4}
For any $i=1,...,n$, we have $P_i(t_1,...,t_s)=\sum_{\alpha=1} ^s L_i(f_{\alpha})^p (t_{\alpha} ^p -t_{\alpha})$.
\end{lem}
\dem For any index $i$, it is easy to check by induction on $N>0$ that the expression:
$$
\left(T_{i,\sum_{\alpha=1} ^s t_{\alpha} L_i(f_{\alpha})}\right) ^{N-1}
\left(\sum_{\alpha=1} ^s t_{\alpha} L_i(f_{\alpha})\right)
$$
is a polynomial in $A[t_1,...,t_s]$ of degree $N$ in $t_1,...,t_s$, whose leading term is equal to:
$$
\left(\sum_{\alpha=1} ^s t_{\alpha} L_i(f_{\alpha})\right)^N
$$
In particular, the polynomial $P_i$ has degree $p$ in $t_1,...,t_s$, and its leading term is given by:
$$
\left(\sum_{\alpha=1} ^s t_{\alpha} L_i(f_{\alpha})\right)^p=\sum_{\alpha=1} ^s t_{\alpha} ^p L_i(f_{\alpha})^p
$$
Moreover $P_i$ vanishes at all points of $\FP^s$ by lemma \ref{solution2}. By lemma \ref{solution3}, there exist
some elements $a_{i,1},...,a_{i,s}$ of $A$ such that:
$$
P_i(t_1,...,t_s)=a_{i,1} (t_1 ^p - t_1) + ...+ a_{i,s} (t_{s} ^p - t_{s})
$$
Since the leading term of $P$ is equal to $\sum_{\alpha} L_i(f_{\alpha})^p t_{\alpha} ^p$, we have $L_i(f_{\alpha})^p=a_{i,\alpha}$
for any $\alpha=1,...,s$ and the result follows.
\qed

\subsection{Proof of Theorem \ref{solutions}}

Let $f_1,...,f_s$ be some algebraic solutions of $\{D_1,...,D_n\}$ such that $L(f_1),...,L(f_s)$ form a basis
of $Vect_k(I)$. We are going to show that $L(f_1),...,L(f_s)$ form a basis of $I$ over $\mathbb{F}_p$.
Since $\GAF$ is isomorphic to $I$, Theorem \ref{solutions} will follow. By construction, these elements are linearly
independent over $k$ (hence over $\mathbb{F}_p$). So it suffices to prove that they span $I$ over $\mathbb{F}_p$.
Let $x$ be any element of $I$. By definition of $I$, there exists an algebraic solution $f$ of $\{D_1,...,D_n\}$ such that
$x=L(f)$. Since the $L(f_j)$ form a basis of the $k$-vector space spanned by $I$, there exist some elements $z_1,...,z_s$
of $k$ such that $L(f)=z_1 L(f_1)+...+z_s L(f_s)$. We only need to check that
$z_1,...,z_s$ belong to $\FP$. By Lemmas \ref{solution2} and \ref{solution4}, we have for any $i=1,...,n$: 
$$
P_i(z_1,...,z_s)= \sum_{\alpha=1} ^s L_i(f_{\alpha})^p (z_{\alpha} ^p -z_{\alpha})=0
$$
But $k$ is algebraically closed. So for any $\alpha=1,...,s$, there exists an element $y_{\alpha}$ of $k$ such that
$y_{\alpha} ^p=z_{\alpha} ^p -z_{\alpha}$. Since $char(k)=p$, we get for any $i=1,...,n$:
$$
\sum_{\alpha=1} ^s y_{\alpha} ^p L_i(f_{\alpha})^p=\left(\sum_{\alpha=1} ^s y_{\alpha} L_i(f_{\alpha}) \right)^p=0
$$
Since $A$ is an integral $k$-algebra, it has no nilpotent elements and $\sum_{\alpha=1} ^s y_{\alpha} L_i(f_{\alpha}) =0$
for any index $i$. But this is equivalent to writing:
$$
\sum_{\alpha=1} ^s y_{\alpha} L(f_{\alpha})=0
$$
Since $L(f_1),...,L(f_s)$ are linearly independent over $k$, we get $y_1=...=y_{s}=0$. Therefore
$z_{\alpha} ^p -z_{\alpha}=0$ for any $\alpha=1,...,s$, every $z_{\alpha}$ belongs to $\FP$ and
Theorem \ref{solutions} is proved.

\subsection{Proof of Corollary \ref{estimate}}

Let $A$ be a $k$-algebra of finite type over an algebraically closed field $k$ of characteristic $p$. Let $\FF$ be a foliation on
$Spec(A)$, and denote by $D_1,...,D_n$ a set of generators of $\FF$ as an $A$-module. We follow the notations of the previous subsection
with $\DD=\FF$. Let $\deg$ be a good degree function on $A$
(see the Introduction). For any integer $n$, let ${\cal{L}}(n)$ be the $k$-vector space of the elements
$f$ of $A$ such that $\deg(f)\leq n$. Let $L$ stand for the morphism defined at the beginning of this section
for $D_1,...,D_n$. Then we have the inclusion:
$$
I=L(\Pi(A,\FF))\subseteq {\cal{L}}(\deg(D_1))\times ...\times {\cal{L}}(\deg(D_n))
$$
by definition of the degree of a $k$-derivation. In particular, we find:
$$
\dim_k Vect_k(I)\leq l(\deg(D_1))+...+l(\deg(D_n))
$$ 
Let $f_1,...,f_s$ be some algebraic solutions of $\FF$ such that $L(f_1),...,L(f_s)$ form a basis over $k$ of $Vect_k(I)$. By the arguments
of the previous
subsection, $L(f_1),...,L(f_s)$ form also a basis of $I$ over $\FP$. In particular we have:
$$
\dim_{\FP} I=\dim_k Vect_k(I)
$$
Since $L$ is an isomorphism between $\Pi(A,\FF)$ and $I$, we obtain the inequality:
$$
\dim_{\FP} \Pi(A,\FF) \leq l(\deg(D_1))+...+l(\deg(D_n))
$$
which is exactly the assertion of Corollary \ref{estimate}.

\section{Proof of Corollary \ref{picard}}

Let $A$ be an integral $k$-algebra of finite type over an algebraically closed field $k$ of characteristic $p>0$. Let $B$ be a subalgebra
of $A$ such that $A^p \subseteq B\subseteq A$. We assume that $A$ and $B$ are normal. Consider the pull-back morphism $\pi$ defined by:
$$
\pi: Pic(B) \longrightarrow Pic(A), \quad M\longmapsto M\otimes_B A
$$
Since $A,B$ are normal, there exists a foliation $\FF_B$ on $Spec(A)$ such that $B=\cap_{D\in \FF_B} \ker D$ (see \cite{M-P}).
Denote by $G$ {\em the group of algebraic solutions of $\FF_B$ modulo the subgroup spanned by the invertible elements of $A$
and the first integrals of $\FF_B$}. Note that every invertible element of $A$ is an algebraic solution.
By construction, we have a surjective morphism:
$$
F: \Pi(A,\FF_B) \longrightarrow G
$$
Since $\Pi(A,\FF_B)$ is finite dimensional over $\FP$ by Theorem \ref{solutions}, $G$ is also finite dimensional. We
are going to construct an injective group morphism:
$$
\theta: \ker\pi \longrightarrow G
$$
This will imply that $\ker\pi$ is finite dimensional over $\FP$, and give the estimate:
$$
dim_{\FP} \ker\pi\leq dim_{\FP} G \leq dim_{\FP} \Pi(A,\FF_B)
$$
which is exactly the result given by Corollary \ref{picard}. \\

The construction proceeds as follows. Let $M'$ be a finitely generated, locally free $B$-submodule of $K(B)$. If $M'\otimes A\simeq A$, then
$A.M'\simeq A$ and there
exists an element $f'$ of $K(A)$ such
that $A.M'=A\{f'\}$. Since $K(B)$ is annihilated by $\FF_B$, $D(M')=0$ for any $D\in \FF_B$ and the $A$-module $A.M'$ is stable by $\FF_B$.
In particular,
for any element $D$ of $\FF_B$, $D(f')$ belongs to $A\{f'\}$ and $f'$ is an algebraic solution of $\FF_B$. Note that $f'$ is uniquely
determined up to
multiplication by an element of $A^*$. If ${\cal{E}}$ denotes the set of finitely generated, locally free $B$-submodules $M'$ of $K(B)$ such that
$M'\otimes A\simeq A$, then we have a well-defined correspondence:
$$
\Theta: {\cal{E}}\longrightarrow G, \quad M' \longmapsto [f']
$$
Let $M''$ be another finitely generated, locally free $B$-submodule of $K(B)$. If $M''\simeq M'$, then the isomorphism from $M'$ to $M''$ is
induced by the multiplication by an element $g$ of $K(B)^*$, i.e. $M''=g.M'$. Since $A.M''=A\{f''\}$, we have:
$$
A.M''=A\{f''\}=g.A.M'=A\{gf'\}
$$
and $f''/gf'$ is an invertible element of $A$. So the class $[f']$ of $f'$ in $G$ only depends on the isomorphism class of $M'$. Since every
invertible sheaf $M$ on $Spec(B)$ can be represented by a finitely generated, locally free $B$-submodule $M'$ of $K(B)$, $\Theta$ induces
a map:
$$
\theta: \ker\pi  \longrightarrow G, \quad [M'] \longmapsto [f']
$$
\begin{lem}
The map $\theta$ is a group morphism.
\end{lem}
\dem Let $M',M''$ be two finitely generated, locally free $B$-submodules of $K(B)$ such that $M'\otimes A\simeq A$ and $M''\otimes A\simeq A$. Let
$M'M''$ be the $B$-submodule of $K(B)$ spanned by all products $x'x''$, where $x'\in M'$ and $x''\in M''$. Since $M'\otimes A$ and $M''\otimes A$
are trivial, there exist some elements $f',f''$ of $K(A)$ such that $A.M'=A\{f'\}$ and $A.M''=A\{f''\}$. But then we have:
$$
A.M'M''=(A.M')(A.M'')=A\{f'\}.A\{f''\}=A\{f'f''\}
$$
and the $A$-module $A.M'M''$ is spanned by $f'f''$. Since the module $M'M''$ represents the invertible sheaf $M' \otimes M''$ in $Pic(B)$, we obtain:
$$
\theta(M'\otimes M'')=\theta(M'M'')=[f'f'']=[f'][f'']=\theta(M')\theta(M'')
$$
\qed

\begin{lem} \label{injection}
The morphism $\theta$ is injective.
\end{lem}
\dem Let $M'$ be a finitely generated, locally free $B$-submodule of $K(B)$. Assume that $M'\otimes A\simeq A$ and that $[f']=0$. Then
there exists a rational first integral $g$ of $\FF_B$, and an invertible element $h$ of $A$, such that $A.M'=A\{gh\}=A\{g\}$. Since
$D(g)=0$ for any $D\in \FF_B$, $g$ belongs to $K(B)$. First we show that $M'\subseteq B\{g\}$. Let $x$ be any element of $M'$. Since $x$
belongs to $A.M'$, there exists an element $a$ of $A$ such that $x=a.g$. Since $D(x)=D(g)=0$ for any $D\in \FF_B$, $D(a)=0$ for any
$D\in \FF_B$ and $a$ belongs to $B$. In particular, $x$ belongs to $B\{g\}$ and we have:
$$
M'\subseteq B\{g\}
$$ 
Second we show the equality $M'=B\{g\}$. Since $M'$ and $B\{g\}$ are finite modules over a noetherian ring, it suffices to
prove that ${M'}_{\cal{M}}=B_{\cal{M}}\{g\}$ for any maximal ideal
${\cal{M}}$ of $B$. Let ${\cal{M}}$ be a maximal ideal of $B$. Since $A^p \subseteq B\subseteq A$, $A$ is integral over $B$. By the Cohen going-up
Theorem (see \cite{Ei}), there exists a maximal ideal ${\cal{M}}_A$ of $A$ such that ${\cal{M}}_A\cap B={\cal{M}}$. Since $M'$ is locally free of rank 1, there exists
an element $x$ of $M'$ such that ${M'}_{\cal{M}}=B_{\cal{M}}\{x\}$. But then we find:
$$
A_{{\cal{M}}_A}M'=A_{{\cal{M}}_A}{M'}_{\cal{M}}=A_{{\cal{M}}_A}B_{\cal{M}}\{x\}=A_{{\cal{M}}_A}\{x\}
$$
On the other hand, we have:
$$
A_{{\cal{M}}_A}M'=(A.M')_{{\cal{M}}_A}=A_{{\cal{M}}_A}\{g\}
$$
In particular, the fraction $a=x/g$ is an invertible element of $A_{{\cal{M}}_A}$. Since $x$ and $g$ belong to $K(B)$, $a$ and
$a^{-1}$ belong to $K(B)$.
But $A^p \subseteq B\subseteq A$, so $A_{{\cal{M}}_A}$ is a finite $B_{\cal{M}}$-module. In particular, $a$ and $a^{-1}$ are
integral over $B_{\cal{M}}$.
Since $B$ is normal, $B_{\cal{M}}$ is also normal and $a, a^{-1}$ belong to $B_{\cal{M}}$. Therefore $a$ is an invertible
element of $B_{\cal{M}}$. Since
${M'}_{\cal{M}}\subseteq B_{\cal{M}}\{g\}$, we obtain the equality:
$$
{M'}_{\cal{M}}=B_{\cal{M}}\{x\}=B_{\cal{M}}\{g\}
$$
Since this holds for any maximal ideal ${\cal{M}}$, we have $M'=B\{g\}$. In particular, the invertible sheaf $M$ represented by
$M'$ is trivial, and injectivity follows.
\qed

\section{Two examples}

\subsection{Computation of $\Pi(A,\FF)$}

let $A$ be the $k$-algebra $k[x,y]$, where $k$ is algebraically closed of characteristic $p$. Note that $A$ is normal and $A^*=k^*$.
Let $\deg$ be the standard homogeneous degree on $A$. Given an element $t$ of $k$, consider the following $k$-derivation $D_t$ on $A$:
$$
D_t= x\frac{\partial}{\partial x} + ty \frac{\partial}{\partial y}
$$
By an easy compution, we find $(D_t)^p=D_{t^p}$. Then two cases may occur. \\ \\
{\it \underline{First case}: $t$ belongs to $\FP$.} \\ \\
Then the foliation $\FF(D_t)$ is equal to $A\{D_t\}$. Since $D_t$ has degree zero, the dimension of
$\Pi(A, \FF(D_t))$ is bounded by $l(0)=1$ by Corollary \ref{estimate}. Since $D_t(x)=x$, $x$ is an
algebraic solution whose class in $\Pi(A, \FF(D_t))$ is nonzero. Therefore:
$$
\dim_{\FP} \Pi(A,\{D_t \})=\dim_{\FP} \Pi(A, \FF(D_t))=1
$$
{\it \underline{Second case}: $t$ does not belong to $\FP$.} \\ \\
Then $(D_t)^p$ and $D_t$ are $A$-linearly independent,
and the foliation $\FF(D_t)$ contains the $A$-module $A\{x\partial/\partial x, y\partial/\partial y\}$.
Since this latter is $p$-closed, stable by Lie bracket and contains $D_t$, we have $\FF(D_t)=A\{x\partial/\partial x, y\partial/\partial y\}$.
So $\FF(D_t)$ is spanned by 2 derivations of degree zero. By Corollary \ref{estimate}, the dimension of $\Pi(A, \FF(D_t))$ is bounded
by $l(0)+l(0)=2$. Consider now the polynomials $x$ and $y$. Since $D_t(x)/x=1$, $D_t(y)/y=t$ and $t\not\in \FP$, $x$ and $y$
are $\FP$-linearly independent in $\Pi(A, \FF(D_t))$. Therefore:
$$
\dim_{\FP} \Pi(A,\{D_t \})=\dim_{\FP} \Pi(A, \FF(D_t))=2
$$
Note that for all $t$, the Picard group of the kernel $B$ of $\FF(D_t)$ is reduced to zero. Therefore
the estimate given in Theorem \ref{picard} is not the best one in case $t\in \FP^*$. In either cases, we can
notice that the dimension of $\Pi(A, \DD)$ depends not only on the degrees of the elements of $\DD$,
but also on those of the whole foliation $\FF(\DD)$.

\subsection{Computation of a Picard group}

let $A$ be the $k$-algebra $k[x,y]$, where $k$ is algebraically closed of characteristic $2$. Let $B$ be the subalgebra
$k[x^2,y^2,x(1+xy)]$. Then $B$ is isomorphic to $k[u,v,w]/(w^2 -u(1+uv))$. This latter ring is normal because the variety
$Spec(B)$ is a complete intersection which is nonsingular in codimension 1 (see \cite{Ha}). So $B$ is a normal subring
of $A$, which is itself normal. We would like to compute the Picard group of $B$ by means of Corollary \ref{picard}.
To that purpose, consider the following $k$-derivation $D$:
$$
D=x^2 \frac{\partial}{\partial x} + \frac{\partial}{\partial y}
$$
Since $D^2=0$, the foliation $\FF(D)$ is equal to $A.\{D\}$. Moreover $D$ vanishes on $B$, hence $B\subseteq B'=\ker D$. Since
$K(B)\subseteq K(B')\subseteq K(A)$ and that $K(A)/K(B)$ is an extension of degree 2, $K(B')$ is either equal to $K(B)$ or to
$K(A)$. The latter case is impossible because $D$ is a nonzero derivation, and so $K(B')=K(B)$. Since $B$ and $B'$ contain $A^2$,
$B'$ is integral over $A^2$, hence over $B$. As $B'\subseteq K(B)$ and $B$ is normal, we deduce that $B'=B$. In particular, we
have:
$$
\FF_B=A.\{D\}
$$
Now $A$ is a UFD, $A^*=k^*$ and $B$ is normal. By Lemma \ref{injection}, the morphism $\theta$ injects $Pic(B)$ in $\Pi(A,\{D\})$. First we determine
$\Pi(A,\{D\})$. Let $\deg$ be the standard homogeneous degree on $k[x,y]$. Consider the injective morphism:
$$
L:\Pi(A,\{D\})\longrightarrow A, \quad f\longmapsto \frac{D(f)}{f}
$$
Since $D$ has degree 1, $L(f)$ is of the form $a+bx+cy$, where $a,b,c$ belong to $k$. By lemma \ref{solution1},
this polynomial satisfies the equation $D(L(f))=L(f)^2$, which implies:
$$
D(bx+cy)=bx^2 + c=(a+bx+cy)^2
$$
So $a=0$, $b$ belongs to $\mathbb{F}_2$ and $c=0$. In particular, $dim_{\mathbb{F}_2} \Pi(A,\{D\})\leq 1$. Since $D(x)=x^2$, $\Pi(A,\{D\})$
is nonzero of dimension 1, spanned by $x$. Consider now the ideal $M'$ of $B$ defined by $M'=(x^2,x(1+xy))$. The open sets $D(x^2)$ and $D(1+x^2y^2)$
build a covering of $Spec(B)$, and we have:
$$
{M'}_{(1/x^2)}=B_{(1/x^2)} \quad {\rm{and}} \quad {M'}_{(1/1+x^2y^2)}=B_{(1/1+x^2y^2)}\{x\}
$$
Therefore $M'$ defines an invertible sheaf $M$ on $Spec(B)$. Since $A.M'=A\{x\}$ and that $D(x)\not=0$, we find $\theta(M)=[x]\not=0$. In particular
$\theta$ is an isomorphism from $Pic(B)$ to $\Pi(A,\{D\})$, and $Pic(B)={\mathbb{F}_2}$.


\begin{thebibliography}{Jou2}

\bibitem[B-N]{B-N} M.Brunella, M.Nicolau {\it Sur les hypersurfaces solutions des \'equations de Pfaff}, C. R. Acad. Sci. Paris S\'er. I
Math., 329 (1999), $n^o$ 9, 793-795.

\bibitem[Ei]{Ei} D.Eisenbud {\it Commutative Algebra with a view toward Algebraic
Geometry}, Springer Verlag New York (1995).

\bibitem[Ek]{Ek} T.Ekedahl {\it Foliations and inseparable morphisms}, Algebraic Geometry, Bowdoin, 1985 (Brunswick, Maine, 1985), 139-149, Proc. Symp. Pure
Math., 46, Part 2, Amer. Math. Soc., Providence, RI, 1987.

\bibitem[Ha]{Ha} R.Hartshorne {\it Algebraic Geometry}, Graduate Texts in Mathematics $n^o$ 52, Springer Verlag New York Heidelberg (1977).

\bibitem[Jou]{Jou} J-P.Jouanolou {\it Equations de Pfaff alg\'ebriques}, Lect. Notes in Math. 708, Springer Verlag Berlin (1979).

\bibitem[Jou2]{Jou2} J-P.Jouanolou {\it Hypersurfaces solutions d'une \'equation de Pfaff analytique}, Math. Ann. 232 (1978) $n^o$ 3, 239-245.

\bibitem[Mo]{Mo} S.D.Morrison {\it Continuous derivations}, J. Algebra 110, 468-479 (1987).

\bibitem[M-P]{M-P} Y.Miyaoka, T.Peternell {\it Geometry of higher-dimensional varieties}, DMV Seminar 26, Birkh\"auser Verlag, Basel 1997.

\bibitem[Pe]{Pe} J.Pereira {\it Invariant hypersurfaces for positive characteristic vector fields},
J. Pure Appl. Algebra 171 (2002), $n^o$ 2-3, pp. 295-301.

\bibitem[Sh]{Sh} I.Shafarevich {\it Basic Algebraic Geometry}, second edition, Springer Verlag Berlin 1994.
J. Pure Appl. Algebra 171 (2002), $n^o$ 2-3, pp. 295-301.

\end{thebibliography}
\end{document}